\documentclass[aop,MSNbibl,seceqn,dvips]{arximspdf}
\usepackage{mathbh}

\doi{10.1214/12-AOP807} %kopijuoti is PTS
\volume{42}
\issue{2}
\pubyear{2014}
\firstpage{559}
\lastpage{575}

\makeatletter

\newcommand{\cal}{\mathcal}

\newtheorem{Theorem}{Theorem}[section]
\newproclaim{Definition}[Theorem]{Definition}
\newtheorem{Proposition}[Theorem]{Proposition}
\newtheorem{Corollary}[Theorem]{Corollary}
\newproclaim{Remark}[Theorem]{Remark}

\def\e{\mathbb E}
\def\p{\mathbb P}
\def\R{\mathbb R}

\makeatother

\begin{document}
\begin{frontmatter}

\title{Characterization of positively correlated squared Gaussian processes}
\runtitle{Positively correlated squared Gaussian}

\begin{aug}
\author[A]{\fnms{Nathalie} \snm{Eisenbaum}\corref{}\ead[label=e1]{nathalie.eisenbaum@upmc.fr}}
\runauthor{N. Eisenbaum}
\affiliation{CNRS, Universit\'e Pierre et Marie Curie}
\address[A]{Laboratoire de probabilit\'es\\
\quad et mod\`eles al\'eatoires\\
Universit\'e Paris 6\\
Case 188-4, Place Jussieu-75252 Paris cedex 05\\
France\\
\printead{e1}} %adresu isvedimo komanda gale!
\end{aug}

% HISTORY:
\received{\smonth{1} \syear{2012}}
\revised{\smonth{9} \syear{2012}}

% ABSTRACT
%
\begin{abstract}
We solve a conjecture raised by Evans in 1991 on the characterization
of the positively correlated squared Gaussian vectors. We extend this
characterization from squared Gaussian vectors to permanental vectors.
As side results, we obtain several equivalent formulations of the
property of infinite divisibility for squared Gaussian processes.
\end{abstract}

% KEYWORDS

\end{frontmatter}

%s1 #&#
\section{Introduction}

A random vector $(\psi_j)_{1 \leq j \leq n}$ of $\R^n$ is said to be
``associated'' or ``positively correlated'' if for every couple of
increasing functions $F,G$ from $\R^n$ into $\R$ (i.e., $F$ and $G$
are increasing in each variable)
%
%e1.1 #&#
\begin{equation}
\label{poscor} \e\bigl( FG\bigl( (\psi_j)_{1 \leq j \leq n}\bigr)\bigr)
\geq \e\bigl(F\bigl((\psi_j)_{1 \leq j
\leq n}\bigr)\bigr) \e\bigl(G
\bigl((\psi_j)_{1 \leq j \leq n}\bigr)\bigr).
\end{equation}
In 1982, Pitt \cite{Pitt} has shown that a centered Gaussian vector
$\eta= (\eta_i)_{1 \leq i \leq n}$ is ``positively correlated'' iff
the entries of its covariance matrix are all nonnegative, which means
that the Gaussian vector is positively correlated in the usual sense.
To distinguish between the two meanings for positive correlation, we
will keep the writing ``positively correlated,'' in inverted commas, to
refer to the definition (\ref{poscor}).

In 1991, Evans \cite{Ev1} conjectured that given a centered Gaussian
vector $\eta= (\eta_i)_{1 \leq i \leq n}$, the squared centered
Gaussian vector $\eta^2 = (\eta_i^2)_{1 \leq i \leq n}$ is
``positively correlated'' iff there exists a function $\sigma$ from $\{
1 \leq i \leq n\}$ into $\{ -1, 1\}$ such that $(\sigma(i) \eta_i)_{1
\leq i \leq n}$ is positively correlated.

We prove the following:
%
%th1.1 #&#
\begin{Theorem}\label{solution}
A squared centered Gaussian vector is ``positively correlated''
if and only if it is infinitely divisible.
\end{Theorem}

Evans condition for a squared centered Gaussian vector to be
``positively correlated'' is hence necessary but not sufficient.\vadjust{\goodbreak} Indeed, several
necessary and sufficient conditions for a squared centered Gaussian
vector to be infinitely divisible have been established that allow to
see this. The first one was found by
Griffiths \cite{G} in 1983, simplified then by Bapat \cite{B}. This
condition has been translated in terms of Green function of Markov
processes by Eisenbaum and Kaspi \cite{EK}. Another version of this
condition has been established by Vere-Jones \cite{V}. We will use
Vere-Jones characterization of infinitely divisible squared Gaussian
vectors to establish three other equivalent necessary and sufficient
conditions for a squared centered Gaussian process with continuous
covariance to be ``positively correlated.'' One extends the definition
(\ref{poscor}) from vectors to processes by saying that a process is
``positively correlated'' if all its finite-dimensional marginals are
``positively correlated.''

Eisenbaum and Kaspi's characterization stems from the desire to
understand which were the Gaussian processes involved in Dynkin's
isomorphism theorem~\cite{D1}. Here is a brief presentation of the
content of this theorem. Consider a symmetric transient Markov process
$X$ with state space $E$ and $0$-potential density (i.e., Green
function) $(g(x,y), (x,y) \in E\times E)$. The function $g$ is positive
definite. Denote by $(\eta_x)_{x \in E}$ a centered Gaussian process
with covariance $g$, independent of $X$. For $a$, $b$ in $E$ such that
$g(a,b) > 0$, denote by $\tilde{\p}_{ab}$ the probability under which
$X$ starts at $a$ and dies at its last visit to $b$. Besides, $X$
admits a local time process. Denote by $(\tilde{L}^{ab}(x), x \in E)$
the process of the total accumulated local times under $\tilde{\p
}_{ab}$. Then\vspace*{1pt} according to Dynkin's isomorphism theorem, the process
$(\tilde{L}^{ab}(x) + {1\over2}\eta^2_x, x \in E)$ has the same law
as $({1\over2} \eta^2_x, x \in E)$ under the measure ${1\over\e
[\eta_a \eta_b]}\e[\eta_a \eta_b,\cdot]$.\vspace*{3pt}

This identity in law immediately raises two questions: Which are the
centered Gaussian processes with a covariance equal to a Green
function? Which are the centered Gaussian processes $\eta$ such that
the law of $\eta^2$ under $\e[\eta_a \eta_b,\cdot]$ is a positive measure?

An answer to the first question has been given in \cite{EK} (completed
then in \cite{E1}; see (\ref{bibi}) below) under the following form.
Given\vspace*{1pt} a centered Gaussian process $(\eta_x)_{x \in E}$ with a
continuous positive definite covariance $(G(x,y), (x,y) \in\break  E \times E)$,
\textit{$(\eta^2_x)_{x \in E}$ is infinitely divisible if and only if
there exist a real nonnegative
measurable function $d$ on $E$ and a function $g$ on $E^2$ such that}
%
%e1.2 #&#
\begin{equation}
G(x,y) = d(x) g(x,y) d(y)
\end{equation}
\textit{and $g$ is the Green function of a symmetric transient Markov process.}

The corollary below actually provides three alternative formulations to
this answer. One of them is our solution to Evans conjecture for
processes. Another one answers also to the second question. To
introduce the remaining one, we will use the following definition.
%
%de1.2 #&#
\begin{Definition}\label{D1} A random process $(\phi_t)_{t \in E}$ is
said to satisfy \textit{Fortuyin Kasteleyn Ginibre's inequality} (FKG
inequality) if for some reference positive measure $m$, for every
integer $n$, every $t = (t_1,t_2,\ldots, t_n)$ in $E^n$, $(\phi_{t_1},
\phi_{t_2},\break \ldots, \phi_{t_n})$ has a density with respect to $m(dx_1)\cdots m(dx_n)$ product measure on $\R^n$
denoted by $h_t$ such that for every $x$, $y$ in $\R^n$,
\[
h_t(x)h_t(y) \leq h_t(x\wedge y)
h_t(x \vee y),
\]
where $x\wedge y = (x_1 \wedge y_1, x_2 \wedge y_2,\ldots,x_n \wedge y_n)$ and
$x \vee y = (x_1 \vee y_1, x_2 \vee y_2,\ldots,\allowbreak x_n \vee y_n)$.
\end{Definition}
%
%co1.3 #&#
\begin{Corollary}\label{T1} Let $(\eta_x)_{x \in E}$ be a centered
Gaussian process with a continuous positive definite covariance
$(G(x,y), (x,y) \in E\times E)$. The following four properties are equivalent:
\begin{longlist}
\item[(1)] $\eta^2$ is infinitely divisible.

\item[(2)] $\eta^2$ is ``positively correlated.''

\item[(3)] $\eta^2$ satisfies the FKG inequality.

\item[(4)] For every $(a,b)$ in $E^2$, the law of $\eta^2$ under
$\e[ \eta_a \eta_b,\cdot]$ is a positive measure.
\end{longlist}
\end{Corollary}
Once the question of the characterization of ``positively correlated''
squared centered Gaussian processes is solved, one may ask the same
question for shifted Gaussian processes. In particular, given a
centered Gaussian process $ (\eta_x)_{x \in E}$ and a real number $r$,
when is the process $ ((\eta_x + r)^2)_{x \in E}$ ``positively
correlated''?
Thanks to \cite{BW} and \cite{Ev}, we know a sufficient condition for
the realization of that property: the infinite divisibility of $((\eta_x +r)^2)_{x \in E}$. But there is not known characterization, in terms
of the covariance of $\eta$, of that condition for a fixed $r$.
Nevertheless, in \cite{E1}, we have established the following
characterization. Assuming that the set $E$ contains more that two
elements (see Remark \ref{not2}), let $(\eta_x)_{x \in E}$ be a
centered Gaussian process with a continuous covariance
%
%e1.3 #&#
\begin{eqnarray}
\label{bibi} &\mbox{\textit{$\bigl((\eta_x + r)^2
\bigr)_{x \in E}$ is infinitely divisible for every real $r$}},&
\nonumber\\
&\mbox{\textit{if and only if}}&
\\
&\mbox{\textit{the covariance of $\eta$ is the Green function of a transient
Markov process.}}&\nonumber
\end{eqnarray}

This will be used to enunciate another sufficient condition for $((\eta_x + r)^2)_{x \in E}$ to be ``positively correlated'' for every $r$.

The paper is organized as follows. In Section \ref{I} we prove Theorem
\ref{solution}. We then deduce Corollary \ref{T1}. The proofs involve
stochastic comparison of squared centered Gaussian vectors. As a side
result, for a given covariance $G$, we give necessary and sufficient
conditions for the stochastic monotonicity of the family of squared
Gaussian vectors with the resolvents of $G$ for respective
covariance.\vadjust{\goodbreak}

In Section \ref{II} we extend our characterization of ``positively
correlated'' squared Gaussian vectors to permanental vectors. This
extension is legitimated by the fact that a connection, similar to
Dynkin isomorphism theorem, has been established in~\cite{EK}, between
permanental processes and local times of not necessarily symmetric
Markov processes.

In Section \ref{III} we establish an equivalent formulation of (\ref{bibi}).

As it will be shown in Sections \ref{I}, \ref{II} and \ref{III},
many properties of Gaussian processes and, more generally, of
permanental processes, are hence conditioned to the fact that their
kernel is a Green function or not. So it is interesting to mention a
way to generate Green functions. This is done in Section \ref{green}.

%s2 #&#
\section{\texorpdfstring{Proof of Theorem \protect\ref{solution} and Corollary \protect\ref{T1}}
{Proof of Theorem 1.1 and Corollary 1.3}}\label{I}

The proof of Theorem \ref{solution} will show some other equivalent
properties to infinite divisibility for squared Gaussian vectors. To
formulate them, we make use
of the following definitions.
%
%de2.1 #&#
\begin{Definition} A random vector $(\phi_i)_{1 \leq i \leq n}$ of $\R^n$ \textit{stochastically dominates} another random vector $(\psi_i)_{1
\leq i \leq n}$ of $\R^n$ if for any increasing function $F$ from $\R^n$ into $\R$,
\[
\e\bigl[F(\phi_{1},\phi_{2},\ldots, \phi_{n})
\bigr] \geq\e\bigl[F(\psi_{1},\psi_{2},\ldots,
\psi_{n})\bigr].
\]
\end{Definition}
%
%de2.2 #&#
\begin{Definition}\label{strong} Let $(\phi_{t_1},\phi_{t_2},\ldots,
\phi_{t_n})$ and $(\psi_{t_1},\psi_{t_2},\ldots, \psi_{t_n})$ be
two random vectors of $\R^n$, such that there exists a positive
measure $m$ on $\R$, such that their laws both admit respective
densities $h$ and $f$ with respect to $m(dx_1)\cdots m(dx_n)$ product measure on $\R^n$.
If for every $x$, $y$ in $\R^n$,
\[
f(x) h(y) \leq f(x \wedge y) h(x \vee y),
\]
then one says that $(\phi_{t_1},\phi_{t_2},\ldots, \phi_{t_n})$ is
strongly stochastically bigger than (or strongly stochastically
dominates) $(\psi_{t_1},\psi_{t_2},\ldots, \psi_{t_n})$.

One extends this definition to a couple $(\phi,\psi)$ of real-valued
processes indexed by the same set by saying that $\phi$ is strongly
stochastically bigger than $\psi$ when all the finite-dimensional
marginals of $\phi$ and $\psi$ satisfy the above relation.
\end{Definition}

Strong stochastic domination implies usual stochastic domination.
%
%de2.3 #&#
\begin{Definition} Let $C$ be a positive semi-definite matrix. For
$\alpha> 0$, one defines the associated $\alpha$-resolvent matrix as
$C_{\alpha} = (I + \alpha C)^{-1} C$.
\end{Definition}

We have the following corollary of Theorem \ref{solution}.
%
%co2.4 #&#
\begin{Corollary}\label{cor1} Let $\eta= ( \eta_i)_{1 \leq i \leq
n}$ be a centered Gaussian vector with covariance $G$, an $n\times
n$-positive definite matrix. Denote\vadjust{\goodbreak} by $\eta_{\alpha} = (\eta_{\alpha}(i))_{1 \leq i \leq n}$
a centered Gaussian vector with
covariance $G_{\alpha}$. Then the four following points are equivalent:
\begin{longlist}[(iii)]
\item[(i)] $\eta^2$ is infinitely divisible.

\item[(ii)] The family of vectors $(\eta^2_{\alpha})_{\alpha
\geq0}$ is stochastically decreasing as $\alpha$ increases on $\R^+$.

\item[(iii)] The family of vectors $(\eta^2_{\alpha})_{\alpha
\geq0}$ is strongly stochastically decreasing as $\alpha$ increases
on $\R^+$.

\item[(iv)] For every couple $(i,j)$, $1 \leq i,j \leq n$, for every $n\times n$ diagonal matrix $D$,
$(\e [|\eta_{\alpha}(i) \eta_{\alpha}(j)| ])_{\alpha\geq0}$ is
decreasing as $\alpha$ increases on $\R^+$ when $G$ is replaced by~DGD.
\end{longlist}
\end{Corollary}

We adopt the following notation from the paper \cite{PittPerl}. For
$C$ a $n\times n$-positive definite matrix and any measurable function
$F$ on $\R^n$,
$\e_{C}[F(\eta)]$ denotes the expectation with respect to a centered
Gaussian vector $\eta$ with covariance matrix $C$.
\begin{pf*}{Proof of Theorem \ref{solution}}
Thanks to \cite{BW} or \cite{Ev}, we know that if the vector $\eta^2$
is infinitely divisible, then it is ``positively correlated.'' We prove
now the converse.

Assume that $\eta^2$ is ``positively correlated.'' Denote by $G =
(G(i,j))_{1 \leq i,j \leq n}$ its covariance matrix. For every
decreasing function $F, H$ on $\R^n$ we have
\[
\e_G\bigl(FH\bigl(\eta^2\bigr)\bigr) \geq
\e_G\bigl(F\bigl(\eta^2\bigr)\bigr) \e_G
\bigl(H\bigl(\eta^2\bigr)\bigr)
\]
and, in particular, for every $\alpha, \varepsilon> 0$
%
%e2.1 #&#
\begin{equation}
\label{6} \e_G\bigl( e^{- ((\alpha+ \varepsilon)/2) \sum_{i = 1}^n \eta^2_i}\bigr) \geq\e_G
\bigl( e^{- (\alpha/2) \sum_{i = 1}^n \eta^2_i}\bigr) \e_G\bigl( e^{- (\varepsilon/2)\sum_{i = 1}^n \eta^2_i}\bigr).
\end{equation}
Moreover, for any decreasing function $F$ on $\R^n_+$, we have
%
%e2.2 #&#
\begin{eqnarray}
\label{epsilon}
&&
\e_G\bigl(F\bigl( \eta^2\bigr)
e^{- ((\alpha+ \varepsilon)/2) \sum_{i =
1}^n \eta^2_i}\bigr) \nonumber\\[-8pt]\\[-8pt]
&&\qquad\geq\e_G\bigl(F\bigl(\eta^2\bigr)
e^{- (\alpha/2) \sum_{i = 1}^n \eta
^2_i}\bigr) \e_G\bigl( e^{- (\varepsilon/2)\sum_{i = 1}^n
\eta^2_i}\bigr).\nonumber
\end{eqnarray}
We make use now of a remark of Marcus and Rosen (Remark 5.2.4, page 200
in~\cite{MR1}) according to which for all measurable function $K$ on
$\R^n$,
%
%e2.3 #&#
\begin{equation}
\label{LemmaA} \e_{G}\bigl[K(\eta) e^{- {(\alpha/2)} \sum_{i = 1}^n \eta^2_i}\bigr] =
\e_{G_{\alpha}}\bigl[K(\eta)\bigr] \e_{G}\bigl[e^{- {(\alpha/2)} \sum_{i = 1}^n
\eta^2_i}
\bigr].
\end{equation}
In particular, we have
%
%e2.4 #&#
\begin{equation}
\label{bibi1} \e_{G}\bigl[F\bigl(\eta^2\bigr)
e^{- {(\alpha/2)} \sum_{i = 1}^n \eta^2_i}\bigr] = \e_{G_{\alpha}}\bigl[F\bigl(\eta^2\bigr)
\bigr] \e_{G}\bigl[e^{- {(\alpha/2)} \sum_{i =
1}^n \eta^2_i}\bigr].
\end{equation}
We mention that, unlike for (\ref{LemmaA}), one does not need to
assume that $G$ is invertible to obtain (\ref{bibi1}) (for a direct
proof see the proof of Proposition \ref{P1} in Section \ref{II}).

Thanks to (\ref{bibi1}), (\ref{epsilon}) can be rewritten as
\begin{eqnarray*}
&&
\e_{G_{\alpha+ \varepsilon}}\bigl[F\bigl( \eta^2\bigr)\bigr] \e_G
\bigl[e^{- ((\alpha+ \varepsilon)/2) \sum_{i = 1}^n \eta^2_i}\bigr] \\
&&\qquad\geq\e_{G_{\alpha}}\bigl[F\bigl(\eta^2
\bigr)\bigr] \e_G\bigl[ e^{- (\alpha/2) \sum
_{i = 1}^n \eta^2_i}\bigr] \e_G
\bigl( e^{- (\varepsilon/2)\sum_{i =
1}^n \eta^2_i}\bigr).
\end{eqnarray*}
Consequently, for every increasing function $F$, we obtain, thanks to
(\ref{6}),
\[
\e_{G_{\alpha+ \varepsilon}}\bigl[F\bigl( \eta^2\bigr)\bigr] \leq
\e_{G_{\alpha
}}\bigl[F\bigl(\eta^2\bigr)\bigr] {\e_G[ e^{- (\alpha/2) \sum_{i = 1}^n \eta
^2_i}] \e_G( e^{- (\varepsilon/2)\sum_{i = 1}^n \eta
^2_i})\over\e_G[e^{- ((\alpha+ \varepsilon)/2) \sum_{i =
1}^n \eta^2_i}]
}.
\]
Thanks to (\ref{6}), we finally obtain for every increasing,
nonnegative function $F$ on~$\R^n_+$,
%
%e2.5 #&#
\begin{equation}
\label{inc} \e_{G_{\alpha+ \varepsilon}}\bigl[F\bigl( \eta^2\bigr)\bigr] \leq
\e_{G_{\alpha
}}\bigl[F\bigl(\eta^2\bigr)\bigr].
\end{equation}
Because of the restriction on the sign of $F$, the above inequality
does not mean stochastic domination but will be sufficient for our
purpose. Indeed, for a fixed $\alpha> 0$, note that
\[
G_{\alpha+ \varepsilon} = ( I + \varepsilon G_{\alpha})^{-1}
G_{\alpha}.
\]
Set $f_{\alpha} (\varepsilon) = \e_{G_{\alpha+ \varepsilon}}[F(
\eta^2)]$, and note that $f_{\alpha}$ is decreasing at $0$.

Besides, we set $C_{ij}(G_{\alpha}) = G_{\alpha}(i,j)$. We also
define a function $\cal F$ on the set of covariance matrices by setting
\[
{\cal F}(C) = \e_C\bigl[F\bigl(\eta^2\bigr)\bigr].
\]
In \cite{PittPerl}, the derivatives of functions of the form $\e_C[H(\eta)]$ with respect to the entries of the matrix are computed.
The authors work with a $C^2(\R^n)$-function $H$ which together with
its first and second derivatives satisfy a $O(|x|^N)$ growth condition
at $\infty$, for some finite $N$. For $F$ measurable function on
$\R^n_+$ such that the function $H$ defined by $H(x_1,\ldots, x_n) =
F(x_1^2,\ldots, x_n^2)$ satisfies this condition, one easily obtains
for $i\not= j$,
%
%e2.6 #&#
\begin{equation}
\label{1} {\partial{\cal F} \over\partial C_{ij}} (C) = 4 \e_C\biggl[
\eta_i \eta_j \,{\partial^2 F\over\partial x_i \,\partial x_j}\bigl(
\eta^2\bigr)\biggr]
\end{equation}
and
%
%e2.7 #&#
\begin{equation}
\label{2} {\partial{\cal F} \over\partial C_{ii}} (C) = 2 \e_C\biggl[
\eta_i^2 \,{\partial^2 F\over\partial x^2_i }\bigl(\eta^2
\bigr)\biggr] + \e_C\biggl[{\partial F\over
\partial x_i }\bigl(
\eta^2\bigr)\biggr].
\end{equation}
For $\varepsilon$ small enough, we have
$G_{\alpha+ \varepsilon} = \sum_{k = 0}^{\infty} (-1)^k \varepsilon^k (G_{\alpha})^{k+1}$, hence,
%
%e2.8 #&#
\begin{equation}
\label{3} C_{ij}(G_{\alpha+ \varepsilon}) = \sum_{k = 0}^{\infty}(-1)^k
\varepsilon^k (G_{\alpha})^{k+1}(i,j),
\end{equation}
which is a derivable function of $\varepsilon$ at $0$. We obtain
%
%e2.9 #&#
\begin{equation}
f'_{\alpha}(\varepsilon) = \sum_{1 \leq i\leq j \leq n}
{\partial
{\cal F} \over\partial C_{ij}} (G_{\alpha+ \varepsilon}) \,{\partial
C_{ij} \over\partial\alpha}
(G_{\alpha+ \varepsilon}),
\end{equation}
which, thanks to (\ref{1}), (\ref{2}) and (\ref{3}), leads to
%
%e2.10 #&#
\begin{eqnarray}
\label{derive} f'_{\alpha}(0) &=& -4 \sum
_{1 \leq i < j \leq n} \e_{G_{\alpha
}}\biggl[\eta_i
\eta_j \,{\partial^2 F\over\partial x_i \,\partial x_j}\bigl(\eta^2\bigr)\biggr]
(G_{\alpha})^2(i,j)
\nonumber\\[-8pt]\\[-8pt]
&&{} - \sum_{i = 1}^n \e_{G_{\alpha}}
\biggl[ 2\eta_i^2 \,{\partial^2 F\over
\partial x^2_i }\bigl(
\eta^2\bigr) + {\partial F\over\partial x_i }\bigl(\eta^2\bigr)
\biggr] (G_{\alpha})^2(i,i)\nonumber
\end{eqnarray}
[we mention that $(G_{\alpha})^2(i,j)$ is not $(G_{\alpha}(i,j))^2$].

We choose now to take $F(x) = \sqrt{x_ix_j}$, with $i \not= j$. We
first check that (\ref{1}) and (\ref{2}) still hold. Indeed, the
formulas computed in \cite{PittPerl} are still available for $H(x) =
|x_i x_j|$. For this choice (\ref{derive}) gives
%
%e2.11 #&#
\begin{equation}
\label{in} \e_{G_{\alpha}}\bigl[\operatorname{sgn}( \eta_i\eta_j)
\bigr] (G_{\alpha})^2(i,j) \geq0.
\end{equation}
Note that for every $(\lambda_k)_{1 \leq k \leq n}$ in $\R^n$, the
vector $(\lambda_k^2 \eta^2_k)_{1 \leq k \leq n}$ is also
``positively correlated.'' Consequently, setting $\lambda=
\operatorname{Diag}((\lambda_k)_{1 \leq k \leq n})$, one can replace $G_{\alpha}$ by
$\lambda G_{\alpha} \lambda$ in (\ref{in}) to
obtain
\[
\operatorname{sgn}(\lambda_i \lambda_j) \e_{G_{\alpha}}\bigl( \operatorname{sgn}(
\eta_i \eta_j)\bigr) \lambda_i
\lambda_j \sum_{k = 1}^n
G_{\alpha}(i,k) \lambda_k^2 G_{\alpha}(k,j)
\geq0,
\]
which is equivalent to
\[
\sum_{k = 1}^n \lambda_k^2
\e_{G_{\alpha}}\bigl(\operatorname{sgn}(\eta_i \eta_j)\bigr)
G_{\alpha}(i,k) G_{\alpha}(k,j)\geq0.
\]
Since this is true for every $\lambda$, we have
%
%e2.12 #&#
\begin{equation}
\label{4}\qquad \e_{G_{\alpha}}\bigl(\operatorname{sgn}(\eta_i \eta_j)
\bigr) G_{\alpha}(i,k) G_{\alpha
}(k,j) \geq0 \qquad\mbox{for every }
i,j,k \mbox{ with } i \not=j.
\end{equation}
We choose to take $k = i$ and obtain
%
%e2.13 #&#
\begin{equation}
\label{couple} G_{\alpha}(i,j) \e_{G_{\alpha}}\bigl[ \operatorname{sgn}(
\eta_i \eta_j)\bigr] \geq0,
\end{equation}
which together with (\ref{4}) leads to
%
%e2.14 #&#
\begin{equation}
\label{5} G_{\alpha}(j,i) G_{\alpha}(i,k) G_{\alpha}(k,j) \geq0
\qquad\mbox{for every } i,j,k \mbox{ with } i \not=j.
\end{equation}
We show now that this condition implies that there exists $\sigma_{\alpha}$ from $\{1,2,\ldots,n \}$ into $\{-1, 1 \}$ such that for every $i,j$,
%
%e2.15 #&#
\begin{equation}
\label{Evans} \sigma_{\alpha}(i) G_{\alpha}(i,j)
\sigma_{\alpha}(j) \geq0.
\end{equation}
We do it by recurrence on the size of the matrix $G_{\alpha}$. Assume
that our claim is true at rank $n$ and suppose that $G_{\alpha}$ is a
$(n+1)\times(n+1)$-covariance matrix. We just need to define $\sigma_{\alpha}(n+1)$. For every $j, k$ in $\{1,2,\ldots,n\}$, we have \mbox{$\sigma_{\alpha}(j)\sigma_{\alpha}(k) G_{\alpha}(j,k) \geq0$}.
Since $G_{\alpha}(n+1,j) G_{\alpha}(j,k) G_{\alpha}(k,n+1) \geq
0$, we obtain
$\sigma_{\alpha}(j)\sigma_{\alpha}(k) G_{\alpha}(n+1,j) G_{\alpha
}(n+1,k) \geq0$. Consequently,
$\sigma_{\alpha}(j) G_{\alpha}(n+1,j)$ has a constant sign
independent of $j$, $1 \leq j \leq n$, that we denote by $\sigma_{\alpha}(n+1)$. This implies immediately that
$\sigma_{\alpha}(j) G_{\alpha}(n+1,j) \sigma_{\alpha}(n+1) \geq0$.

We then easily check that our claim holds for $n = 3$.

For a real positive number $\beta$ and a $m\times m$-matrix $M =
(M_{i,j})_{1 \leq i,j \leq m}$, the quantity $\mathrm{per}_{\beta}(M)$
is defined as follows: $\mathrm{per}_{\beta}(M) = \sum_{\tau\in{\cal
S}_m} \beta^{\nu(\tau)} \prod_{i = 1}^m M_{i,\tau(i)}$
where ${\cal S}_m$ is the set of the permutations on $\{1,2,\ldots,m\}$,
and $\nu(\tau)$ is the signature of $\tau$.

For every integer $m$, every $k_1, k_2,\ldots, k_m$ in $\{1,2,\ldots, n\}$
and every $\beta> 0$, we hence have
\begin{eqnarray*}
\mathrm{per}_{\beta} \bigl( \bigl(G_{\alpha}(k_i,k_j)
\bigr)_{1\leq i,j \leq m}\bigr) &=&\sum_{\tau\in{\cal S}_m}
\beta^{\nu(\tau)} \prod_{i = 1}^n
G_{\alpha}(k_i,k_{\tau(i)})
\\
&=& \sum_{\tau\in{\cal S}_m} \beta^{\nu(\tau)} \prod
_{i = 1}^n \sigma_{\alpha}(k_i)
\sigma_{\alpha}(k_{\tau(i)})G_{\alpha
}(k_i,k_{\tau(i)})
\geq0,
\end{eqnarray*}
which is a sufficient condition for $\eta^2$ to be infinitely
divisible thanks to the Vere-Jones criteria \cite{V} (this criteria is
recalled at the beginning of Section~\ref{II}).
\end{pf*}
\begin{pf*}{Proof of Corollary \ref{T1}}
One can easily notice that (1) is equivalent to~(3). Indeed, according
to Bapat \cite{B}, a centered Gaussian vector $(\eta_i)_{1 \leq i
\leq n}$ with nonsingular covariance matrix $G$ is such that $(\eta^2_i)_{1 \leq i \leq n}$ is infinitely divisible iff there exists a
signature matrix $\sigma$ [a diagonal matrix such that $\sigma(i,i) =
-1$ or~$1$] such that $\sigma G^{-1} \sigma$ is a $M$-matrix (i.e.,
its off-diagonal entries are nonpositive). Thanks to \cite{KR1}, we
know that this is also a necessary and sufficient condition for $(\eta^2_i)_{1 \leq i \leq n}$ to satisfy the Fortuyin--Kasteleyn--Ginibre's
inequality. Note that there is no need of the continuity of the
covariance to then conclude on the equivalence between (1) and (3).

Thanks to Theorem \ref{solution}, we hence immediately have the
equivalence of (1), (2) and (3). Note that we did not have to use the
well-known fact that (3) implies (2)~\cite{FKG}.

Under the assumption of continuity of $G$, we know thanks to \cite{EK}
that (1) is realized iff for every $x$, $y$ in $E$, $G(x,y) = d(x)
g(x,y) d(y)$, with $d$ a nonnegative measurable function on $E$ and
$g$ the Green function of some transient Markov process. Denote by
$(\tilde{\eta}_x, x \in E)$ a centered Gaussian process with covariance
$g$. Thanks to Dynkin's isomorphism theorem, we know that for every
$a$, $b$ in $E$, the law of $(\tilde{\eta}^2_x)_{x \in E}$ under
$\e[\tilde{\eta}_a \tilde{\eta}_b,\cdot]$ is a positive measure.
Since $(\eta_x)_{x \in E} \stackrel{\mathrm{(law)}}{=} (d(x)
\tilde{\eta}_x)_{x \in E}$, we see that (4) is realized.\vadjust{\goodbreak}

To see that (4) implies (1), note first that for every $x$, $y$ in $E$,
$G(x,y) \geq0$. Denote by $\mathbb G$
the matrix $(G(x_i,x_j))_{1
\leq i,j \leq n}$ for $x_1,x_2,\ldots,x_n$ in $E$ and denote by
${\mathbb G}_{\alpha}$ the $\alpha$-resolvent matrix associated to
$\mathbb G$. We note that, thanks to (\ref{LemmaA}), for every $\alpha
>0$, we have for every $a$ and $b$ in $\{x_1,\ldots, x_n\}$,
\begin{eqnarray*}
\e_{\mathbb G}\bigl[\eta_a\eta_b e^{- {(\alpha/2)} \sum_{i = 1}^n
\eta^2_{x_i}}
\bigr] &=& \e_{{\mathbb G}_{\alpha}}[\eta_a\eta_b]
\e_{\mathbb G}\bigl[e^{- {(\alpha/2)} \sum_{i = 1}^n \eta^2_i}\bigr] \\
&=& {\mathbb G}_{\alpha}(a,b)
\e_{\mathbb G}\bigl[e^{- {(\alpha/2)} \sum
_{i = 1}^n \eta^2_i}\bigr].
\end{eqnarray*}
Hence, for every $\alpha> 0$, and every $a$ and $b$, ${\mathbb
G}_{\alpha}(a,b) \geq0$, which, according to Vere-Jones (Proposition
4.5 in \cite{V}, recalled at the beginning of Section~\ref{II}), is a
sufficient condition for $(\eta_{x_i}^2, 1 \leq i \leq n)$ to be
infinitely divisible. Since this is true for every $x_1,\ldots,x_n$, we
conclude that $\eta^2$ is infinitely divisible.
\end{pf*}
\begin{pf*}{Proof of Corollary \ref{cor1}} We start by noting that the
density $f_{\alpha}$ of $\eta^2_{\alpha}$ with respect to the
Lebesgue measure is connected to the density $f_0$ of $\eta^2$.
Indeed, thanks to (\ref{LemmaA}), we have for a.e. $x =
(x_1,x_2,\ldots,x_n)$ in $\R^n_+$,
%
%e2.16 #&#
\begin{equation}
\label{density} f_{\alpha}(x) = {e^{- {(\alpha/2)} \sum_{i = 1}^n x_i }\over\e
[\operatorname{exp}\{- {(\alpha/2)} \sum_{i = 1}^n \eta^2_i\}]}
f_0(x).
\end{equation}
Assume now that (i) is satisfied. Thanks to Corollary \ref{T1}, this
implies that $\eta^2$ satisfies the FKG inequality. Thanks to (\ref
{density}), one obtains for $\alpha< \beta$ and every $x,y$ in $\R^n_+$,
\[
f_{\alpha}(x) f_{\beta}(y) \leq f_{\alpha}(x \vee y)
f_{\beta
}(x\wedge y),
\]
which leads to (iii).

Now (iii) implies (ii) and (ii) implies (iv). The proof of Theorem \ref
{solution} shows that (iv) implies (i).
\end{pf*}

%s3 #&#
\section{The nonsymmetric case}\label{II}

A real-valued positive vector $(\psi_i, 1 \leq i \leq n)$ is a
permanental vector if its Laplace transform satisfies
for every
$(\alpha_1, \alpha_2,\ldots,\break \alpha_n)$ in~$\R^n_+$,
%
%e3.1 #&#
\begin{equation}
\label{perm} \e\Biggl[\operatorname{exp}\Biggl\{-{1\over2} \sum
_{i = 1}^n \alpha_i
\psi_{i}\Biggr\}\Biggr] = | I + \alpha G |^{-1/\beta},
\end{equation}
where $I$ is the $n\times n$-identity matrix, $\alpha$ is the diagonal
matrix $\operatorname{Diag}((\alpha_i)_{1 \leq i \leq n})$,
$G = (G(i,j))_{1\leq i,j \leq n}$ and $\beta$ is a fixed positive number.

Such a vector $(\psi_i, 1 \leq i \leq n)$ is a permanental vector
with kernel
$(G(i,j),\break  1 \leq i,j \leq n)$ and index $\beta$.

Permanental vectors represent a natural extension of squared centered
Gaussian vectors. Indeed, for $\beta= 2$ and $G$ covariance matrix,
(\ref{perm}) is the Laplace transform of a squared centered Gaussian vector.

Thanks to Vere-Jones (Proposition 4.5 in \cite{V}), we know that there
exists a nonnegative random vector with Laplace transform given by
(\ref{perm}) if and only if:

\begin{longlist}[(II)]
\item[(I)]
All the real eigenvalues of $G$ are nonnegative.

\item[(II)] For every $\alpha> 0$, set $G_{\alpha} = (I + \alpha G)^{-1} G$,
then $G_{\alpha}$ is $\beta$-positive.
\end{longlist}

A $n\times n$-matrix $M = (M(i,j))_{1\leq i,j \leq n}$ is said to be
$\beta$-positive if
for every integer~$m$, every $k_1, k_2,\ldots, k_m$ in $\{1,2,\ldots, n\}$
\[
\mathrm{per}_{\beta} \bigl( \bigl(M(k_i,k_j)
\bigr)_{1\leq i,j \leq m}\bigr) \geq0,
\]
where for any $m\times m$-matrix $A = (A(i,j))_{1 \leq i,j \leq m}$,
the quantity $\mathrm{per}_{\beta}(A)$ is defined as follows: $
\mathrm{per}_{\beta}(A) = \sum_{\tau\in{\cal S}_m} \beta^{\nu(\tau)}
\prod_{i = 1}^m A_{i,\tau(i)}$,
with ${\cal S}_m$ the set of the permutations on $\{1,2,\ldots,m\}$, and
$\nu(\tau)$ the signature of $\tau$.

Obviously, a permanental vector with kernel $G$ is infinitely divisible
if and only if it satisfies the Vere-Jones conditions for every $\beta
> 0$.

Note that the kernel of a permanental vector is not uniquely determined.
We have proved in \cite{EK} that a permanental vector is infinitely
divisible iff it admits as kernel the Green function of some
transient Markov process.
%
%th3.1 #&#
\begin{Theorem}\label{T2} Let $\psi$ be a permanental vector with
index $2$ and kernel $G$. The two following properties are equivalent:
\begin{longlist}[(2)]
\item[(1)] $\psi$ is infinitely divisible.

\item[(2)] $\psi$ is ``positively correlated.''
\end{longlist}
\end{Theorem}
To prove Theorem \ref{T2}, we need the following preliminary
proposition, that will be established at the end of this section.
%
%pr3.2 #&#
\begin{Proposition}\label{P1} For $\beta> 0$, let $M$ be a $n\times
n$ matrix such that there exists a random nonnegative vector $\psi=
(\psi(1),\psi(2),\ldots,\psi(n))$ with Laplace transform
\[
\e\bigl(e^{-({1/2}) \sum_{i = 1}^n x_i\psi(i)}\bigr) = |I + x M|^{-1/\beta}
\]
for every $(x_1,x_2,\ldots,x_n)$ in $R^n_+$. Set for every $\alpha\geq
0$, $M_{\alpha} = M(I + \alpha M)^{-1}$.

There exists a nonnegative random vector $\psi_{\alpha} = (\psi_{\alpha}(1),\psi_{\alpha}(2),\ldots,\psi_{\alpha}(n))$ with Laplace transform
\[
\e\bigl(e^{-{(1/2)} \sum_{i = 1}^n x_i\psi_{\alpha}(i)}\bigr) = |I + x M_{\alpha}|^{-1/\beta}.
\]
The law of $\psi_{\alpha}$ is absolutely continuous with respect to
the law of $\psi$. Moreover, for every bounded measurable functional
$F$ on $\R^n_+$, we have
\[
\e\bigl[F(\psi_{\alpha})\bigr] = \e\biggl[{ \operatorname{exp}\{- {(\alpha/2)} \sum_{i = 1}^n
\psi(i) \}\over\e[\operatorname{exp}\{- {(\alpha/2)} \sum_{i = 1}^n \psi(i) \}
]} F(\psi)
\biggr].
\]
\end{Proposition}
\begin{pf*}{Proof of Theorem \ref{T2}} Let $G$ be a $n\times n$-matrix such
that there exists a permanental vector with index $2$ and kernel $G$.
For any measurable function $F$ on~$\R^n_+$,
$\e_{G}[F(\psi)]$ denotes the expectation with respect to a
permanental vector $\psi$ with covariance matrix $G$ and index $2$.

We already know, thanks to \cite{BW} or \cite{Ev}, that (1) implies
(2). We show that (2) implies (1). Assume that $\psi$ is ``positively
correlated.''
Thanks to Proposition \ref{P1}, for every measurable function $F$ on
$\R^n_+$,
%
%e3.2 #&#
\begin{equation}
\label{absolute} \e_{G} \bigl[ F(\psi) e^{-{(\alpha/2)} \sum_{i = 1}^n \psi_i}\bigr] =
\e_{G_{\alpha}} \bigl[ F(\psi)\bigr]\e_G\bigl[ e^{-{(\alpha/2)} \sum_{i = 1}^n
\psi_i}
\bigr].
\end{equation}
Similarly as in the proof of Theorem \ref{solution}, one hence obtains
that for every nonnegative increasing function $F$ on $\R^n_+$,
%
%e3.3 #&#
\begin{equation}
\label{21} \e_{G_{\alpha}}\bigl[F(\psi)\bigr] \leq\e_{G}\bigl[F(
\psi)\bigr].
\end{equation}
Now we use the fact noticed in \cite{V} that for every $i \not= j$,
$G_{ij} G_{ji} \geq0$.
Remark that for every permanental vector $(\psi_i, \psi_j)$ with
index $2$ and kernel
the $2\times2$-matrix $C$, we have for every function $F$,
%
%e3.4 #&#
\begin{equation}
\label{laplace} \e_C\bigl[F(\psi_i,\psi_j)
\bigr] = \e_{\overline{C}}\bigl[F\bigl(\eta^2_i,
\eta^2_j\bigr)\bigr]
\end{equation}
with the covariance matrix $\overline{C}$ defined by
$\overline{C}_{ii} = C_{ii}$, $\overline{C}_{jj} = C_{jj}$ and
$\overline{C}_{ij} = \sqrt{C_{ij}C_{ji}}$.

Indeed, to prove (\ref{laplace}), one just compares the respective
Laplace transform of the two random couples and checks that for every
$2\times2$-diagonal matrix $x$ with nonnegative entries,
\[
| I + x C| = |I + x \overline{C}|.
\]

Choosing $F(x) = \sqrt{x_ix_j}$ on $\R^n_+$, we obtain, thanks to
(\ref{21}),
\[
\e_{G_{\alpha}}[\sqrt{\psi_i\psi_j}] \leq
\e_{G}[\sqrt{\psi_i\psi_j}],
\]
which together with (\ref{laplace}) leads to
\[
\e_{\overline{G}_{\alpha}} \Bigl[\sqrt{\eta^2_i
\eta^2_j}\Bigr] \leq\e_{\overline{G}} \Bigl[\sqrt
{\eta^2_i\eta^2_j}
\Bigr],
\]
where $\overline{G}_{\alpha}$ is the $2\times2$-matrix defined by
$\overline{G}_{\alpha} (i,i) = G_{\alpha}(i,i)$, $\overline
{G}_{\alpha}(j,j) = G_{\alpha}(j,j)$ and $\overline{G}_{\alpha}(i,j) =
\sqrt{G_{\alpha}(i,j) G_{\alpha}(j,i)}$.

Setting $f(\alpha) = \e_{\overline{G}_{\alpha}}
[\sqrt{\eta^2_i\eta^2_j}]$, we know that $f$ is decreasing at $0$.
Using\break  the same arguments as in the proof of Theorem \ref{solution}, for
$\alpha $ small enough, we have $f'(\alpha) = - 4
\e_{\overline{G}_{\alpha}}[ \operatorname{sgn}(\eta_i\eta_j)]
\,{\partial\overline{G}_{\alpha}(i,j)\over\partial\alpha}$, with
${\partial\overline{G}_{\alpha}(i,j)\over\partial\alpha} =
{1\over2} (G_{\alpha}(i,j) G_{\alpha}(j,\break i))^{-1/2}\{ G_{\alpha
}(i,j)\* G'_{\alpha}(j,i) + G_{\alpha}(j,i) G'_{\alpha}(i,j)\}$.

Hence, we obtain
\[
f'(0) = -{1\over\overline{G}(i,j)} \e_{\overline{G}}\bigl[\operatorname{sgn}(
\eta_i\eta_j)\bigr] \bigl\{ G^2(i,j) G(j,i)
+ G^2(j,i) G(i,j)\bigr\}.
\]
Consequently, we must have
\[
\e_{\overline{G}}\bigl[\operatorname{sgn}(\eta_i\eta_j)\bigr] \bigl
\{ G^2(i,j) G(j,i) + G^2(j,i) G(i,j)\bigr\} \geq0.
\]
Note that since the couple $(\eta^2_i, \eta^2_j)$ is always
infinitely divisible, we have,
using (\ref{couple}), $ \e_{\overline{G}}[\operatorname{sgn}(\eta_i\eta_j)] \geq
0$. Hence, we have
%
%e3.5 #&#
\begin{equation}
\label{ineq} G^2(i,j) G(j,i) + G^2(j,i) G(i,j) \geq0.
\end{equation}
Remark that for every $ (\lambda_1, \lambda_2,\ldots, \lambda_n)$ in
$\R^n_+$, the permanental vector
$(\lambda_1 \psi_1,\allowbreak \lambda_2 \psi_2,\ldots,\lambda_n \psi_n)$ is
also ``positively correlated.'' Since
\[
\e\Biggl[\operatorname{exp}\Biggl\{-{1\over2} \sum
_{i = 1}^n \alpha_i \lambda_i
\psi_{i}\Biggr\}\Biggr] = | I + \alpha\lambda G |^{-1/2},
\]
$(\lambda_1 \psi_1, \lambda_2 \psi_2,\ldots,\lambda_n \psi_n)$ admits
$\lambda G$ for the kernel.
In particular, $\lambda G$ satisfies (\ref{ineq}), which gives
\[
\sum_{k=1}^n \lambda_k \bigl
\{G(i,j) G(j,k)G(k,i) + G(j,i) G(i,k) G(k,j) \bigr\} \geq0
\]
and, consequently, we obtain for every $i,j,k$ with $i \not=j$,
\[
\bigl\{G(j,i) G(j,k)G(k,i) + G(j,i) G(i,k) G(k,j) \bigr\} \geq0.
\]
Since $G(i,j)G(j,i) \geq0$, $G(j,k)G(k,j) \geq0$ and $G(i,k)G(k,i)
\geq0$, the two terms
$G(i,j) G(j,k)G(k,i)$ and $G(j,i) G(i,k) G(k,j)$ have the same sign.
Their sum can be nonnegative only if they are both nonnegative. We
have obtained for every $i$, $j$, $k$
\[
G(j,i) G(j,k)G(k,i) \geq0.
\]
By substituting $(\alpha+ \varepsilon)$ to $\alpha$ in (\ref
{absolute}), one obtains similarly for every $\alpha> 0$
\[
G_{\alpha}(j,i) G_{\alpha}(j,k)G_{\alpha}(k,i) \geq0.
\]

We can then develop the same argument as in the proof of Theorem \ref
{solution} from (\ref{5}) until the conclusion that $\psi$ has to be
infinitely divisible.
\end{pf*}
%
%re3.3 #&#
\begin{Remark}
Note that the proof of Theorem \ref{T2} shows
that a permanental vector $(\psi_i)_{1 \leq i \leq n}$ is infinitely
divisible if and only if for every $i$, $j$, $ 1 \leq i,j \leq n$,
for every $n\times n$ nonnegative diagonal matrix $\lambda$, $\e_{(\lambda G)_{\alpha}}
[\sqrt{ \psi_i \psi_j}]$ is a decreasing function of
$\alpha$ on $\R^+$.
\end{Remark}
\begin{pf*}{Proof of Proposition \ref{P1}} We note that
\[
I + x M_{\alpha} = I + xM(I + \alpha M)^{-1} = \bigl(I + (x+
\alpha)M\bigr) (I + \alpha M)^{-1},
\]
where $x + \alpha$ means $(x_1 + \alpha, x_2 + \alpha,\ldots, x_n +
\alpha)$. Taking the determinant of each part of this equation and
then the power $(-1/ \beta)$ gives
\[
| I + x M_{\alpha}|^{-1/\beta} = \e\bigl( X e^{-{(1/2)} \sum_{i =
1}^n x_i \psi(i)}\bigr),
\]
where $X$ is the positive random variable with expectation $1$ defined by
\[
X = \operatorname{exp}\Biggl\{- {\alpha\over2} \sum_{i = 1}^n
\psi(i) \Biggr\}\bigg/ \e\Biggl[\operatorname{exp}\Biggl\{- {\alpha\over2} \sum
_{i = 1}^n \psi(i) \Biggr\}\Biggr].
\]
Hence, $\psi_{\alpha}$ exists and has the law of $\psi$ under $\e(
X,\cdot)$.
\end{pf*}

%s4 #&#
\section{The shifted case}\label{III}

Given a centered Gaussian process $(\eta_x)_{x\in E}$ and a real
number $r$, we write $(\eta+r)^2$ for $((\eta_x + r)^2)_{x \in E}$.
Thanks to \cite{BW} and \cite{Ev}, we know that
%
%e4.1 #&#
\begin{eqnarray}
\label{obvious} &\mbox{\textit{If $(\eta+ r)^2$ is infinitely
divisible for every real $r$},}&\nonumber\\
&\mbox{\textit{then}}&\\
&\mbox{\textit{$(\eta+r)^2$ is} ``\textit{positively correlated}'' \textit{for every
real $r$.}}&\nonumber
\end{eqnarray}

The following theorem gives another sufficient condition for $(\eta+
r)^2$ to be ``positively correlated.''
It can also be seen as an alternative characterization of Gaussian
processes with a covariance equal to the Green function of a Markov
process. We assume that $E$ contains more than two elements.
%
%th4.1 #&#
\begin{Theorem}\label{shifted}
Let $(\eta_x)_{x \in E}$ be a centered Gaussian process with a
continuous positive definite covariance. The following properties are
equivalent:
\begin{longlist}[(2)]
\item[(1)] The covariance of $\eta$ is the Green function of a
transient Markov process.

\item[(2)] The family of processes $((\eta+ r)^2, r \geq0)$ is
strongly stochastically increasing as $r$ increases on $\R^+$.
\end{longlist}
\end{Theorem}
The definition of a strong stochastic comparison is given at the
beginning of Section \ref{I} (Definition \ref{strong}).
\begin{pf*}{Proof of Theorem \ref{shifted}}
$(1) \Longrightarrow(2)$: Assuming (1), we know that for every
positive integer $n$ and every $(x_i)_{1 \leq i \leq n}$ in $E^n$, the
covariance matrix $G$ of the vector $(\eta_{x_i})_{1 \leq i \leq n}$
is the inverse of a diagonally dominant $M$-matrix (see \cite{EK}),
that is, setting $G^{-1} = M$, all the entries of $G$ are nonnegative,
all the off-diagonal entries of $M$ are nonpositive, and for every $k$,
$\sum_{i = 1}^n M_{ki} \geq0$. The fact that $G^{-1}$ is an
$M$-matrix implies that for every $\beta= (\beta_i)_{1 \leq i \leq
n}$ and $\alpha= (\alpha_i)_{1 \leq i \leq n}$ in $\R^n_+$, such
that $\alpha_i \geq\beta_i$, we have, using a result of Fang and Hu
(Theorem 2.3 in \cite{FH}),
\[
\bigl(\bigl(\eta_{x_i} + (G \alpha)_i\bigr)^2
\bigr)_{1 \leq i \leq n} \mbox{ strongly stochastically dominates } \bigl(\bigl(
\eta_{x_i} + (G \beta)_i\bigr)^2
\bigr)_{1 \leq i \leq n}.
\]
Since $G^{-1}$ is diagonally dominant, we know that the vector $G^{-1}
{\mathbh1}$, where $ {\mathbh1}$ is the vector $(1,1,\ldots, 1)^t$ of
$\R^n_+$, belongs to $\R^n_+$. Hence, we can
choose to take $\alpha= r M{\mathbh1}$ and $\beta= r' M{\mathbh1}$,
with $r \geq r'$, to obtain
\[
\bigl((\eta_{x_i} +r)^2\bigr)_{1 \leq i \leq n} \mbox{
strongly stochastically dominates } \bigl(\bigl(\eta_{x_i} +
r'\bigr)^2\bigr)_{1 \leq i \leq n}.
\]
By definition, this means that the sequence of processes $((\eta+ r)^2,
r > 0)$ increases with $r$ with respect to the strong stochastic order.

$(2) \Longrightarrow(1)$: Conversely,\vspace*{1pt} for $r>0$ fixed and $n$
positive integer, denote by $(f_r(x), x \in\R^n_+)$ the density of
the vector $((\eta_{x_i} + r)^2)_{1 \leq i \leq n}$. By assumption for
every $(r,r')$ such that $r > r'$, we have for every $x,y$ in $\R^n_+$,
\[
f_r(x) f_{r'}(y) \leq f_r(x \vee y)
f_{r'}(x\wedge y).
\]
By integrating the above inequality with respect to ${1 \over\sqrt {2\pi}}e^{-r^2/2} \,dr$, one obtains
\[
h(x) f_{r'}(y) \leq h(x \vee y) f_{r'}(x\wedge y),
\]
where $(h(x), x \in\R^n_+)$ is the density of the vector $((\eta_{x_i} + N)^2)_{1 \leq i \leq n}$, with $N $  standard Gaussian
variable independent of $\eta$.

One integrates then this last inequality with respect to $\p( N \in
dr')$ to obtain
\[
h(x) h(y) \leq h(x \vee y) h(x\wedge y),
\]
which means that the vector $((\eta_{x_i} + N)^2, {1 \leq i \leq n})$
satisfies the FKG inequality. Thanks to Theorem \ref{solution}, this
vector is hence infinitely divisible.
Since this is true for every $n $ and every $(x_i)_{1 \leq i \leq n}$,
the process $((\eta_x + N)^2,\break  x \in E)$ is infinitely divisible. We
use now the assumption on the continuity of the covariance of $\eta$
to claim that, thanks to \cite{E1}, this can be so only if the
covariance of $\eta$ is the Green function of a Markov process.
\end{pf*}
%
%re4.2 #&#
\begin{Remark}\label{not2}
The case of Gaussian couples has to be studied as a particular
case. Indeed,
in \cite{E1}, we have shown that, given a centered Gaussian couple
$(\eta_x, \eta_y)$, the couple $((\eta_x +r)^2, (\eta_y+r)^2)$ is
infinitely divisible for every $r$, if and only if
\[
\e(\eta_x\eta_y) \geq0 \quad\mbox{and}\quad \e(
\eta_x\eta_y) \leq\e \bigl(\eta^2_x
\bigr) \e\bigl(\eta_y^2\bigr).
\]
But one can use the two-dimensional case to show that the converse of
(\ref{obvious}) is false. Indeed, consider a centered Gaussian couple
$(\eta_x, \eta_y)$ with covariance matrix $ \bigl({1 \atop
\rho}\enskip{\rho \atop 1} \bigr)$ such that $|\rho| < 1$. Then
according to Corollary 3.1 of Fang and Hu~\cite{FH}, for every $r$
$((\eta_x +r)^2, (\eta_y +r)^2)$ satisfies the FKG inequality. In
particular, \mbox{$((\eta_x +r)^2, (\eta_y +r)^2)$} is ``positively
correlated'' for every $r$. But choosing $\rho< 0$, we see that
$((\eta_x +r)^2, (\eta_y +r)^2)$ cannot be infinitely divisible for
every $r$.
\end{Remark}

%s5 #&#
\section{A stability property for Green functions}\label{green}

%th5.1 #&#
\begin{Theorem}\label{stab} Let $(g(x,y), (x, y )\in E\times E)$ be
the Green function of a transient Markov process. Assume $g$ is
continuous, then for every $\beta\geq1$,
$(g^{\beta}(x,y), (x, y )\in E\times E)$ is also the Green function of
a transient Markov process.
\end{Theorem}
In the case $E$ is finite, the above fact has already been established
by Dellacherie et al. \cite{Del}. To establish the general case, we
first show the following characterization of Green functions, which is
an extension of a result on symmetric Green functions (see Theorems 1.2
and 1.3 in \cite{E1}).
%
%th5.2 #&#
\begin{Theorem}\label{charac} Let $G$ be a continuous function on
$E\times E$. The three following points are equivalent:
\begin{longlist}[(iii)]
\item[(i)] $G$ is the Green function of some Markov process.
\item[(ii)] For every positive real $c$, $G + c$ is the kernel
of an infinitely divisible permanental process.
\item[(iii)] $G + 1$ is the kernel of an infinitely divisible
permanental process.
\end{longlist}
\end{Theorem}
\begin{pf}%{Proof of Theorem \ref{charac}}
We follow the proof of Theorem
1.2 and Theorem 1.3 in \cite{E1}. We insist only on the arguments that
are specific to the nonsymmetric case.

(i) $\Rightarrow$ (ii): Making use of the arguments developed in \cite
{E1}, there exists a recurrent Markov process $X$ such that $G$
represents the $0$ potential densities of $X$ killed at its first
hitting time of $a$, a point outside $E$. We set then $G(a,a) = 0 =
G(a,x) = G(x,a) $ for every $x$ in $E$. We use then an isomorphism
theorem for recurrent Markov processes (Corollary 3.5 in \cite{EK2})
to claim that for every $c > 0$, there exists a permanental process
$(\psi_x, x \in E\cup\{a\})$ with kernel $G + c$ and index $2$,
satisfying for every $r>0$
%
%e5.1 #&#
\begin{equation}
\label{iso} \bigl(\bigl(\tfrac{1}{2}\psi_x, x \in E\cup
\{a \} \bigr)| \psi_a = r\bigr) \stackrel {\mathrm{(law)}} {=}
\bigl(\tfrac{1}{2}\phi_x + L^x_{\tau_r},
x \in E\cup\{ a \}\bigr),
\end{equation}
where $(\phi_x, x \in E \cup\{a \})$ is a permanental process with
kernel $G$ and index $2$ independent of $X$, and $( L^x_{\tau_r}, x
\in E\cup\{a \})$ is the local time process of $X$ starting at $a$, at
time $\tau_r = \inf\{ s \geq0\dvtx  L^a_s > r\}$.

Since $G$ is a Green function, the process $\phi$ is infinitely
divisible (see \cite{EK2}). Besides, one easily checks that
$L_{\tau_r}$ is infinitely divisible. Actually, $(L_{\tau_r})_{r >0}$
is a L\'evy process and for every
$\alpha= (\alpha_i)_{1 \leq i \leq n}$ in $\R^n_+$ and $(x_i)_{1
\leq i \leq n}$ in $(E \cup\{a \})^n$, we have
%
%e5.2 #&#
\begin{equation}
\label{Laplace} \e\Biggl(\operatorname{exp}\Biggl\{- \sum
_{i=1}^n \alpha_i L^{x_i}_{\tau_r}
\Biggr\}\Biggr) = e^{-
r F(G, \alpha)},
\end{equation}
where $F(G, \alpha)$ is a nonnegative constant.

Hence, for every $r>0$, $(\psi| \psi_a = r)$ is also infinitely
divisible. But (ii) requires the infinite divisibility of $\psi$. We
hence integrate (\ref{iso}) with respect to the law of $\psi_a$ to
obtain, thanks to (\ref{Laplace}),
\[
\e\Biggl(\operatorname{exp}\Biggl\{- {1\over2}\sum
_{i=1}^n \alpha_i \psi_{x_i}
\Biggr\}\Biggr) = \e\Biggl(\operatorname{exp}\Biggl\{- {1\over2}\sum
_{i=1}^n \alpha_i
\phi_{x_i}\Biggr\}\Biggr) \e \bigl( e^{- F(G, \alpha) \psi_a}\bigr).
\]
Now, $\psi_a$ has the law of a squared Gaussian variable and is hence
infinitely divisible. Consequently, for every positive $\delta$, there
exists a nonnegative variable $Y_{\delta}$ that we can choose
independent of $X$, such that
\[
\bigl(\e\bigl( e^{- F(G, \alpha) \psi_a}\bigr)\bigr)^{\delta} = \e\bigl(
e^{- F(G, \alpha)
Y_{\delta}}\bigr).
\]
We hence obtain
\[
\e\Biggl(\!\operatorname{exp}\Biggl\{\!- {1\over2}\sum_{i=1}^n\!\alpha_i \psi_{x_i}
\!\Biggr\}\! \Biggr)^{\delta} = \e\Biggl(\!\operatorname{exp}\Biggl\{\!-
{1\over2}\sum_{i=1}^n\!
\alpha_i \phi_{x_i}\!\Biggr\}\!\Biggr)^{\delta} \e\Biggl(\!
\operatorname{exp}\Biggl\{\!- \sum_{i=1}^n\!
\alpha_i L^{x_i}_{\tau_{Y_{\delta}}}\!\Biggr\}\!\Biggr),
\]
which shows the infinite divisibility of $\psi$.

To prove that (iii) implies (i), we can directly use the argument given
in~\cite{E1}, since symmetry is not required there. And, finally, (ii)
obviously implies~(iii).
\end{pf}

The following equivalence will help us to show Theorem \ref{stab}.
%
%th5.3 #&#
\begin{Theorem}\label{concl} Let $G$ be a continuous function on
$E\times E$. Then $G$ is the Green function of a Markov process if and
only if for every finite subset $F$ of $E$ the restriction of $G$ to
$F\times F$ is the Green function of a Markov process.
\end{Theorem}
\begin{pf}%{Proof of Theorem \ref{concl}}
One has to establish it only in the
nonsymmetric case (in the symmetric case it is a consequence of \cite
{EK} and \cite{E1}). The
direct way is known. Conversely, assume that for every finite set $F$,
$G_{|F\times F}$ is a Green function, then thanks to Theorem \ref
{charac}, $(G + 1)_{|F\times F}$ is the kernel of an infinitely
divisible permanental processes with index $2$. Hence, there exists a
permanental process $(\psi_x, x \in E)$ with kernel $(G(x,y) + 1,
(x,y) \in E\times E)$ and index $2$. Thanks to Theorem \ref{charac},
all the finite-dimensional marginals of $\psi$ are infinitely
divisible. Consequently, $\psi$ is infinitely divisible.
This implies, thanks to Theorem \ref{charac}, that $( G(x,y) +
1,(x,y), (x, y )\in E\times E)$ is the Green function of a Markov
process.
\end{pf}
\begin{pf*}{Proof of Theorem \ref{stab}} Thanks to \cite{Del}, we know
that for every finite subset $F$ of $E$, $(g^{\beta}(x,y), (x, y )\in
F\times F)$ is the Green function of a Markov process.
This implies, thanks to Theorem \ref{concl}, that $(g^{\beta}(x,y),
(x, y )\in E\times E)$ is the Green function of a Markov process.
\end{pf*}


%suskaldyti doi

% imsref loaded by lrinkeviciute, 2013-02-20 08:28:11

\printaddresses


\begin{thebibliography}{18}
% BibTex style file: ims.bst, 2013-01-28
% Default style options (sort=0,type=number).
% Used options (sort=1,type=number).

%b1 ###
\bibitem{B}
\begin{barticle}[mr]
\bauthor{\bsnm{Bapat},~\bfnm{R.~B.}\binits{R.~B.}}
(\byear{1989}).
\btitle{Infinite divisibility of multivariate gamma distributions and
  {$M$}-matrices}.
\bjournal{Sankhy\=a Ser. A}
\bvolume{51}
\bpages{73--78}.
\bid{issn={0581-572X}, mr={1065560}}
\bptok{imsref}%
\end{barticle}
\endbibitem

%b2 ###
\bibitem{BW}
\begin{bincollection}[mr]
\bauthor{\bsnm{Burton},~\bfnm{Robert~M.}\binits{R.~M.}} \AND
  \bauthor{\bsnm{Waymire},~\bfnm{Ed}\binits{E.}}
(\byear{1986}).
\btitle{The central limit problem for infinitely divisible random measures}.
In \bbooktitle{Dependence in Probability and Statistics ({O}berwolfach, 1985)}.
\bseries{Progr. Probab. Statist.}
\bvolume{11}
\bpages{383--395}.
\bpublisher{Birkh\"auser}, \blocation{Boston, MA}.
\bid{mr={0899999}}
\bptok{imsref}%
\end{bincollection}
\endbibitem

%b3 ###
\bibitem{Del}
\begin{barticle}[mr]
\bauthor{\bsnm{Dellacherie},~\bfnm{Claude}\binits{C.}},
  \bauthor{\bsnm{Martinez},~\bfnm{Servet}\binits{S.}} \AND
  \bauthor{\bsnm{San~Martin},~\bfnm{Jaime}\binits{J.}}
(\byear{2009}).
\btitle{Hadamard functions of inverse {$M$}-matrices}.
\bjournal{SIAM J. Matrix Anal. Appl.}
\bvolume{31}
\bpages{289--315}.
\bid{doi={10.1137/060651082}, issn={0895-4798}, mr={2496420}}
\bptok{imsref}%
\end{barticle}
\endbibitem

%b4 ###
\bibitem{D1}
\begin{bincollection}[mr]
\bauthor{\bsnm{Dynkin},~\bfnm{E.~B.}\binits{E.~B.}}
(\byear{1984}).
\btitle{Local times and quantum fields}.
In \bbooktitle{Seminar on Stochastic Processes, 1983 ({G}ainesville, {F}la.,
  1983)}.
\bseries{Progr. Probab. Statist.}
\bvolume{7}
\bpages{69--83}.
\bpublisher{Birkh\"auser}, \blocation{Boston, MA}.
\bid{mr={0902412}}
\bptnote{check year}%
\bptok{imsref}%
\end{bincollection}
\endbibitem

%b5 ###
\bibitem{E1}
\begin{barticle}[mr]
\bauthor{\bsnm{Eisenbaum},~\bfnm{Nathalie}\binits{N.}}
(\byear{2005}).
\btitle{A connection between {G}aussian processes and {M}arkov processes}.
\bjournal{Electron. J. Probab.}
\bvolume{10}
\bpages{202--215 (electronic)}.
\bid{doi={10.1214/EJP.v10-238}, issn={1083-6489}, mr={2120243}}
\bptok{imsref}%
\end{barticle}
\endbibitem

%b6 ###
\bibitem{EK}
\begin{barticle}[mr]
\bauthor{\bsnm{Eisenbaum},~\bfnm{Nathalie}\binits{N.}} \AND
  \bauthor{\bsnm{Kaspi},~\bfnm{Haya}\binits{H.}}
(\byear{2006}).
\btitle{A characterization of the infinitely divisible squared {G}aussian
  processes}.
\bjournal{Ann. Probab.}
\bvolume{34}
\bpages{728--742}.
\bid{doi={10.1214/009117905000000684}, issn={0091-1798}, mr={2223956}}
\bptok{imsref}%
\end{barticle}
\endbibitem

%b7 ###
\bibitem{EK2}
\begin{barticle}[mr]
\bauthor{\bsnm{Eisenbaum},~\bfnm{Nathalie}\binits{N.}} \AND
  \bauthor{\bsnm{Kaspi},~\bfnm{Haya}\binits{H.}}
(\byear{2009}).
\btitle{On permanental processes}.
\bjournal{Stochastic Process. Appl.}
\bvolume{119}
\bpages{1401--1415}.
\bid{doi={10.1016/j.spa.2008.07.003}, issn={0304-4149}, mr={2513113}}
\bptok{imsref}%
\end{barticle}
\endbibitem

%b8 ###
\bibitem{Ev}
\begin{barticle}[mr]
\bauthor{\bsnm{Evans},~\bfnm{Steven~N.}\binits{S.~N.}}
(\byear{1990}).
\btitle{Association and random measures}.
\bjournal{Probab. Theory Related Fields}
\bvolume{86}
\bpages{1--19}.
\bid{doi={10.1007/BF01207510}, issn={0178-8051}, mr={1061945}}
\bptok{imsref}%
\end{barticle}
\endbibitem

%b9 ###
\bibitem{Ev1}
\begin{barticle}[mr]
\bauthor{\bsnm{Evans},~\bfnm{Steven~N.}\binits{S.~N.}}
(\byear{1991}).
\btitle{Association and infinite divisibility for the {W}ishart distribution
  and its diagonal marginals}.
\bjournal{J. Multivariate Anal.}
\bvolume{36}
\bpages{199--203}.
\bid{doi={10.1016/0047-259X(91)90057-9}, issn={0047-259X}, mr={1096666}}
\bptok{imsref}%
\end{barticle}
\endbibitem

%b10 ###
\bibitem{FH}
\begin{barticle}[mr]
\bauthor{\bsnm{Fang},~\bfnm{Zhaoben}\binits{Z.}} \AND
  \bauthor{\bsnm{Hu},~\bfnm{Taizhong}\binits{T.}}
(\byear{1997}).
\btitle{Developments on {${\rm MTP}\sb 2$} properties of absolute value
  multinormal variables with nonzero means}.
\bjournal{Acta Math. Appl. Sin. Engl. Ser.}
\bvolume{13}
\bpages{376--384}.
\bid{doi={10.1007/BF02009546}, issn={0168-9673}, mr={1489839}}
\bptok{imsref}%
\end{barticle}
\endbibitem

%b11 ###
\bibitem{FKG}
\begin{barticle}[mr]
\bauthor{\bsnm{Fortuin},~\bfnm{C.~M.}\binits{C.~M.}},
  \bauthor{\bsnm{Kasteleyn},~\bfnm{P.~W.}\binits{P.~W.}} \AND
  \bauthor{\bsnm{Ginibre},~\bfnm{J.}\binits{J.}}
(\byear{1971}).
\btitle{Correlation inequalities on some partially ordered sets}.
\bjournal{Comm. Math. Phys.}
\bvolume{22}
\bpages{89--103}.
\bid{issn={0010-3616}, mr={0309498}}
\bptok{imsref}%
\end{barticle}
\endbibitem

%b12 ###
\bibitem{G}
\begin{barticle}[mr]
\bauthor{\bsnm{Griffiths},~\bfnm{R.~C.}\binits{R.~C.}}
(\byear{1984}).
\btitle{Characterization of infinitely divisible multivariate gamma
  distributions}.
\bjournal{J. Multivariate Anal.}
\bvolume{15}
\bpages{13--20}.
\bid{doi={10.1016/0047-259X(84)90064-2}, issn={0047-259X}, mr={0755813}}
\bptok{imsref}%
\end{barticle}
\endbibitem

%b13 ###
\bibitem{PittPerl}
\begin{barticle}[auto:STB|2013/01/29|08:09:18]
\bauthor{\bsnm{Joag-Dev},~\bfnm{K.}\binits{K.}},
  \bauthor{\bsnm{Perlman},~\bfnm{M.~D.}\binits{M.~D.}} \AND
  \bauthor{\bsnm{Pitt},~\bfnm{L.~D.}\binits{L.~D.}}
(\byear{1983}).
\btitle{Association of normal random variables and Slepian's inequality}.
\bjournal{Ann. Probab.}
\bvolume{11}
\bpages{451--455}.
\bptok{imsref}%
\end{barticle}
\endbibitem

%b14 ###
\bibitem{KR1}
\begin{barticle}[mr]
\bauthor{\bsnm{Karlin},~\bfnm{Samuel}\binits{S.}} \AND
  \bauthor{\bsnm{Rinott},~\bfnm{Yosef}\binits{Y.}}
(\byear{1981}).
\btitle{Total positivity properties of absolute value multinormal variables
  with applications to confidence interval estimates and related probabilistic
  inequalities}.
\bjournal{Ann. Statist.}
\bvolume{9}
\bpages{1035--1049}.
\bid{issn={0090-5364}, mr={0628759}}
\bptok{imsref}%
\end{barticle}
\endbibitem

%b15 ###
\bibitem{MR1}
\begin{bbook}[mr]
\bauthor{\bsnm{Marcus},~\bfnm{Michael~B.}\binits{M.~B.}} \AND
  \bauthor{\bsnm{Rosen},~\bfnm{Jay}\binits{J.}}
(\byear{2006}).
\btitle{Markov Processes, {G}aussian Processes, and Local Times}.
\bseries{Cambridge Studies in Advanced Mathematics}
\bvolume{100}.
\bpublisher{Cambridge Univ. Press}, \blocation{Cambridge}.
\bid{doi={10.1017/CBO9780511617997}, mr={2250510}}
\bptok{imsref}%
\end{bbook}
\endbibitem

%b16 ###
\bibitem{Pitt}
\begin{barticle}[mr]
\bauthor{\bsnm{Pitt},~\bfnm{Loren~D.}\binits{L.~D.}}
(\byear{1982}).
\btitle{Positively correlated normal variables are associated}.
\bjournal{Ann. Probab.}
\bvolume{10}
\bpages{496--499}.
\bid{issn={0091-1798}, mr={0665603}}
\bptok{imsref}%
\end{barticle}
\endbibitem

%b17 ###
\bibitem{V}
\begin{barticle}[mr]
\bauthor{\bsnm{Vere-Jones},~\bfnm{D.}\binits{D.}}
(\byear{1997}).
\btitle{Alpha-permanents and their applications to multivariate gamma, negative
  binomial and ordinary binomial distributions}.
\bjournal{New Zealand J. Math.}
\bvolume{26}
\bpages{125--149}.
\bid{issn={1171-6096}, mr={1450811}}
\bptok{imsref}%
\end{barticle}
\endbibitem

\end{thebibliography}
\end{document}